\documentclass[english,11pt]{article}

\usepackage[english]{babel}
\usepackage{amsmath,amsthm,amssymb,enumerate}
\usepackage[latin1]{inputenc}
\usepackage[T1]{fontenc}

\usepackage{ulem}
\usepackage{psfig,labelfig}
\usepackage{graphicx}
\usepackage[dvips]{epsfig}
\usepackage{lineno}
\usepackage{color}

\newtheorem{thm}{Theorem}[section]
\newtheorem{cor}[thm]{Corollary}
\newtheorem{lem}[thm]{Lemma}
\newtheorem{con}[thm]{Conjecture}
\newtheorem{cla}[thm]{Claim}

\newcommand{\R}{{\mathbb R}}
\newcommand{\T}{{\mathbb T}}

\newcommand{\Z}{{\mathbb Z}}
\newcommand{\N}{{\mathbb N}}

\newcommand{\mn}{(M,g)}

\DeclareMathOperator{\sy}{sys}
\DeclareMathOperator{\syh}{sys_h}

\DeclareMathOperator{\sr}{SR}
\DeclareMathOperator{\srh}{SR_h}

\DeclareMathOperator{\hyp}{hyp}
\DeclareMathOperator{\dist}{dist}

\DeclareMathOperator{\vol}{vol}

\title{Construction of surfaces with large systolic ratio}

\author{Hugo Akrout and Bjoern Muetzel}

\begin{document}
\maketitle

\begin{abstract}
Let $(M,g)$ be a closed, oriented, Riemannian manifold of dimension $m$. We call a \textit{systole} a shortest non-contractible loop in $\mn$ and denote by $\sy(M,g)$ its length. Let $\sr(M,g)=\frac{{\sy\mn}^m}{\vol \mn}$ be the \textit{systolic ratio} of $\mn$. Denote by $\sr(k)$ the supremum of $\sr(S,g)$ among the surfaces of fixed genus $k \neq 0$. In Section 2 we construct surfaces with large systolic ratio from surfaces with systolic ratio close to the optimal value $\sr(k)$ using cutting and pasting techniques. For all $k_i \geq 1$,  this enables us to prove:
\[
\frac{1}{\sr( k_1 + k_2)} \leq  \frac{1}{\sr(k_1)} +  \frac{1}{\sr(k_2)}.
\]
We furthermore derive the equivalent intersystolic inequality for $\srh(k)$, the supremum of the homological systolic ratio. As a consequence we greatly enlarge the number of genera $k$ for which the bound $\srh(k) \geq \sr(k) \gtrsim \frac{4}{9\pi} \frac{\log(k)^2}{k}$ is valid and show that $\srh(k) \leq \frac{(\log(195k)+8)^2}{\pi(k-1)}$ for all $k \geq 76$. In Section 3 we expand on this idea. There we construct product manifolds with large systolic ratio from lower dimensional manifolds.
\end{abstract}

\section{Introduction}
In the present article we denote by a \textit{manifold} a closed, oriented, Riemannian manifold $(M,g)$ of dimension $m \geq 2$. We denote by $(S,g)$ a Riemannian surface. A \textit{systole} of $(M,g)$ is a shortest non-contractible loop. We denote by $\sy(M,g)$ its length. Normalizing by the volume of $\mn$ we obtain
\[
   \sr(M,g)=\frac{{\sy\mn}^m}{\vol \mn},
\]
the \textit{systolic ratio} of $\mn$, which is invariant under scaling of $\mn$. Let
\[
     \sr(k) = \sup \{ \sr(S,g) \mid (S,g) \text{ Riemannian surface of genus } k \neq 0 \}
\]
be the \textit{optimal systolic ratio in genus} $k$. As its reciprocal value is also quite often used in the literature, we call this value $\sigma(k)$ the \textit{optimal systolic area in genus} $k$, i.e. 
\[
\sigma(k)=\inf \{ \frac{\vol(S,g)}{\sy(S,g)^2} \mid (S,g) \text{ Riemannian surface of genus } k \neq 0 \}.
\]
The exact value of $\sr(k)$ is only known for $k=1$. It was proven by Loewner (see \cite{pu}, p. 71) that $\sr(1)=\frac{2}{\sqrt{3}}$. For large $k$ it is known that (see \cite{ks2})
\begin{equation}
    K_1  \frac{\log(k)^2}{k} \leq \sr(k) \leq  K_2  \frac{\log(k)^2}{k},
\label{eq:srk_bound}
\end{equation}
where $K_1$ and $K_2$ are universal, but unknown constants. The best known upper bound is stated in \cite{ks2}, \textbf{Theorem 2.2}:
\begin{equation}
   \sr(k) \leq \frac{1}{\pi} \frac{\log(k)^2}{k}(1 + o(1)) , \text{ when } k \to \infty.
\label{eq:19sr}  
\end{equation}
It was furthermore shown in \cite{bs} that there exists an infinite sequence of genera $(k_i)_i$, such that 
\begin{equation}
\frac{4}{9\pi} \frac{\log(k_i)^2 - c_0}{k_i}  \leq \sr(k_i),
\label{eq:49sr}
\end{equation}
where $c_0$ is a fixed constant. This result comes from the study of hyperbolic surfaces, i.e. of constant curvature $-1$. More families of hyperbolic surfaces satisfying the above inequality can be found in \cite{ksv1}, \cite{ksv2} and \cite{am}.\\
In the case of a surface $(S,g)$ one can also define the \textit{homological systole}, which is a shortest homologically non-trivial loop in $(S,g)$. This is a shortest non-contractible loop that does not separate $(S,g)$ into two parts. We denote by $\syh(S,g)$ its length and define $\srh(S,g) =\frac{\syh(S,g)^2}{\vol (S,g)}$ as the \textit{homological systolic ratio}. Let

\[
     \srh(k) = \sup \{ \srh(S,g) \mid (S,g) \text{ Riemannian surface of genus } k \neq 0 \}
\]
be the \textit{optimal homological systolic ratio in genus} $k$. We call its reciprocal value $\sigma_h(k)$ the \textit{optimal homological systolic area in genus} $k$. It follows immediately that for any surface $(S,g)$ 
\[
      \sy(S,g) \leq \syh(S,g), \text{ \ \ hence \ \ }  \sr(k) \leq \srh(k).
\]  
Hence $\srh(k)$ has the same lower bound as $\sr(k)$ and it follows from \cite{gr2}, \textbf{Theorem 2.C} that $\srh(k)$  satisfies an upper bound of order $\frac{\log(k)^2}{k}$. In this article we show that $\srh(k)$ is smaller than $\frac{(\log(195k)+8)^2}{\pi (k-1)}$ (see \textbf{Theorem 1.3-3}). An open question is, whether $\sr(\cdot)$ and $\srh(\cdot)$ are monotonically decreasing functions with respect to the genus. Though we can not prove or disprove this result, we can at least show the following intersystolic inequalities:

\begin{thm} Let $\sr(k)$ and $\srh(k)$ be the supremum of the systolic ratio and the homological systolic ratio among all closed, oriented, Riemannian surfaces of genus $k \geq 1$. Let $\sigma(k)$ and $\sigma_h(k)$ be the optimal systolic area and homological systolic area in genus $k \geq 1$. Then for all $k_i \geq 1$
\begin{itemize}
\item[1.]  $\sr(k_1 +k_2) \geq  \left( \frac{1}{\sr(k_1)} + \frac{1}{\sr(k_2)}\right)^{-1}$  \text{ \ \ or equally \ \ }  $\sigma(k_1 +k_2) \leq  \sigma(k_1) + \sigma(k_2)$.
\item[2.] $\srh(k_1 +k_2) \geq  \left( \frac{1}{\srh(k_1)} + \frac{1}{\srh(k_2)}\right)^{-1}$ \text{ \ \ or equally \ \ }   $\sigma_h(k_1 +k_2)\leq  \sigma_h(k_1) + \sigma_h(k_2)$.
\end{itemize}
\label{thm:intermsys}
\end{thm}
In \cite{bb}, p. 159, Babenko and Balacheff provide the equivalent inequality of \textbf{Theorem \ref{thm:intermsys}-1} for connected sums of manifolds of dimension $m \geq 3$. The above inequalities imply that $\sr(k)$ and $\srh(k)$ are at least of order $\frac{1}{k}$. Specializing on metrics $\hyp$ of constant curvature minus one, we obtain the optimal systolic ratio for compact hyperbolic surfaces: 
\[
     \sr(k,\hyp) = \sup \{ \sr(S,\hyp) \mid (S,\hyp) \text{ hyperbolic surface of genus } k \neq 0 \},
\]
and define in an analogous manner the optimal homological systolic ratio for compact hyperbolic surfaces. Again we denote by $\sigma(k,\hyp)$ the inverse. In this case the supremum is attained (see \cite{mu}) and much more is known about the corresponding maximal surfaces than in the general case (see \cite{ak},\cite{am}, \cite{ba},\cite{ge},\cite{sc1}, \cite{sc2}  and \cite{sc3}). Especially 
\[
\sr(k,\hyp) = \srh(k,\hyp) \text{ \ \ (see \cite{pa}, \textbf{Theorem 1.1})}.
\]
This equation enables us to show the first statement of the following theorem:
\begin{thm} Let $\sr(k)$ and $\srh(k)$ be the supremum of the systolic ratio and the homological systolic ratio among all closed, oriented, Riemannian surfaces of genus $k \geq 1$ and $\sr(k,\hyp)$ the supremum of the systolic ratio among compact hyperbolic surfaces of genus $k \geq 2$. Let furthermore $\sigma(k)$, $\sigma_h(k)$ and $\sigma(k,\hyp)$ be the corresponding optimal systolic area. We have:
\begin{itemize}
\item[1.]  $ \frac{\sr(k+1)}{1 -\frac{\sqrt{3}\sr(k+1)}{2}} \geq   \sr(k) \geq  \frac{\sr(k+1,\hyp)}{1 +\frac{\sr(k+1,\hyp)}{\pi}}$ \text{ \ or equally \ }  $ \sigma(k+1) - \frac{\sqrt{3}}{2}  \leq   \sigma(k) \leq  \sigma(k+1,\hyp)+ \frac{1}{\pi}$.
\item[2.]  $ \frac{\srh(k+1)}{1 -\frac{\sqrt{3}\srh(k+1)}{2}} \geq   \srh(k) \geq  \frac{\srh(k+1)}{1 +\frac{\srh(k+1)}{\pi}}$  \text{  \ or equally  \ }  $ \sigma_h(k+1) - \frac{\sqrt{3}}{2}  \leq   \sigma_h(k) \leq  \sigma_h(k+1)+ \frac{1}{\pi}$.
\item[3.]  $ \frac{\sr(k+1)}{\sr(k)}  \leq    1+ \frac{\srh(k+1)}{\pi} $ \text{  \ or equally  \ } $ \frac{\sigma(k+1)}{\sigma(k)}  \geq    \frac{\sigma_h(k+1)}{\sigma_h(k+1) + \frac{1}{\pi}}$.
\end{itemize}
\label{thm:intermsys2}
\end{thm}

\textbf{Theorem \ref{thm:intermsys}} and \textbf{\ref{thm:intermsys2}} are obtained by constructing surfaces with large systolic ratio from surfaces with systolic ratio close to the optimal value $\sr(k)$, $\srh(k)$ or $\sr(k,\hyp)$ using cutting and pasting techniques. As a result, we obtain the following statement:
If $(S,g)$ is a surface of genus $k$, such that $\sr(S,g) \geq \frac{4}{9\pi} \frac{\log(k)^2 - c_0}{k} $, then
\begin{equation*}
\sr(k + j_1) \geq \frac{4}{9\pi} \frac{\log( k + j_1)^2- c_1}{k + j_1} 
\text{ \ and \ }  \sr(j_2\cdot k) \ge \frac{4}{9\pi} \frac{\log( j_2 \cdot k)^2- c_2}{j_2 \cdot k}   \text{ \ \  for all }  j_1, j_2 \ll k.
\end{equation*}
Here the second inequality follows from \textbf{Theorem \ref{thm:intermsys}-1} by induction. This suggests that the bound  $\srh(k) \geq \sr(k) \gtrsim \frac{4}{9\pi} \frac{\log(k)^2}{k}$ is valid for a large number of genera.\\
\\
Furthermore, \textbf{Theorem \ref{thm:intermsys}} and \textbf{\ref{thm:intermsys2}} allow us to provide new lower bounds for $\sr(k)$ for small genera $k$. In \textit{Table~\ref{tab:RS_large_sys}} we give a summary of Riemannian surfaces of genus $1 \leq k \leq 25$ with maximal known systolic ratio. Most of these are constructed from the examples presented in \cite{ca},\cite{ck}, \cite{ksv1}, \cite{sc1} and \cite{sc3} using \textbf{Theorem \ref{thm:intermsys}-1}. As the proof is constructive, the lower bound for $\sr(k)$ is attained in the thus constructed surfaces. The best known upper bounds for $\sr(k)$ in \textit{Table \ref{tab:RS_large_sys}} are due to the following sources (see also \cite{ka}, Chapter 11 for a summary):
\begin{itemize}
\item[-]  genus 2: $\sr(2) \leq \frac{2}{\sqrt{3}}$, \cite{ks1}, \textbf{Theorem 1.3}
\item[-]  genus 3-16: $\sr(k) \leq \frac{4}{3}$ for $k \neq 1$, \cite{gr1}, \textbf{Corollary 5.2.B} 
\item[-]  genus 17-25: 
\begin{equation}
\text{\ \ for all \ } r \in (0,\frac{1}{8}), \frac{\log(2r^2 \sr(k))^2}{4\pi\sr(k)\left(\frac{1}{2}-4r\right)^2(k-1)} \geq 1 , \text{ \ \ \cite{ks2},\cite{ka}, inequality (11.4.1). }
\label{eq:srk_detail}
\end{equation}
\end{itemize}

\begin{table}[htbp]
\begin{center}
\begin{tabular}{|c|c|c|c|c|c|}
\hline
\multicolumn{1}{|c|}{genus $k$} & \multicolumn{1}{|c|}{ \parbox[t]{3.8cm}{surface
(name and/or \\ constructed from)}} &\multicolumn{1}{|c|} {\parbox[t]{2.3cm}{lower bound \\ for $\sr(k)$}}  &\multicolumn{1}{|c|} {\parbox[t]{3.0cm}{upper bound for \\ $\sr(k)$ and $\srh(k)$}} &\multicolumn{1}{|c|} {\parbox[t]{2.8cm}{reference for  \\ the lower bound}} \\ \hline

1 & $\T^2_{hex}$ & 1.15 & 1.15  &   \cite{pu} \\
2 & $R_2$  & 0.80& 1.15      &   \cite{ck}, Fig. 2.1 \\
3 & $R_3$ & 0.66 & 1.33  &  \cite{cal}  \\
4 & $R_4$ & 0.60 & 1.33  &   \cite{ck}, Fig. 2.1  \\
5 & $S_5$ & 0.48 & 1.33 &  \cite{sc3} \\
6 & $I_6$ & 0.42 & 1.33 &  \cite{ca} \\
7 & $H_7$& 0.45 & 1.33 &  \cite{ksv1} \\
8 &$R_8$ (via  $\T^2_{hex},H_7)$ & 0.32 & 1.33  &  Th. 1.1-1 \\
9 & via $R_2,H_7$ & 0.29 & 1.33  &   Th. 1.1-1 \\
10 & via $R_3,H_7$ & 0.27 & 1.33  &   Th. 1.1-1 \\
11 & $I(x|z)$ & 0.28 & 1.33 &  \cite{sc1}\\
12 & via $S_5,H_7$& 0.23& 1.33  &   Th. 1.1-1 \\
13 & via $H_{14}$  & 0.23 & 1.33  & Th. 1.2-1 \\
14 & $H_{14}$ & 0.25 & 1.33 &  \cite{ksv1} \\
15 & via $\T^2_{hex},H_{14}$ & 0.21 & 1.33  &   Th. 1.1-1 \\
16 & via $H_{17}$& 0.26 & 1.33 &    Th. 1.2-1 \\
17 & $H_{17}$ & 0.29 & 1.27 &  \cite{ksv1} \\
18 & via $\T^2_{hex},H_{17}$ & 0.23 & 1.22  &   Th. 1.1-1 \\
19 & via $R_2,H_{17}$ & 0.21 & 1.16  &   Th. 1.1-1 \\
20 & via $R_3,H_{17}$ & 0.20 & 1.12  &   Th. 1.1-1 \\
21 & via $R_4,H_{17}$ & 0.19 & 1.08  &   Th. 1.1-1 \\
22 & via $S_5,H_{17}$ & 0.18 & 1.04  &   Th. 1.1-1 \\
23 & via $I_6,H_{17}$ & 0.17 & 1.00 &   Th. 1.1-1 \\
24 & via $H_7,H_{17}$ & 0.17 & 0.97  &   Th. 1.1-1 \\
25 & via $R_8,H_{17}$ & 0.15 & 0.94  &   Th. 1.1-1 \\  \hline
\end{tabular}
\end{center}
\caption{Upper and lower bounds for $\sr(k)$ and $\srh(k)$ in genus $1 \leq k \leq 25$.}
\label{tab:RS_large_sys}
\end{table}

Revisiting the ideas of  the proof of \cite{gr2}, \textbf{\textbf{Theorem 2.C}},  we also show that 
\begin{thm} Let $\sr(k)$ and $\srh(k)$ be the supremum of the systolic ratio and the homological systolic ratio among all closed, oriented, Riemannian surfaces of genus $k \geq 1$. Then
\begin{itemize}
\item[1.]  $\srh(k) \leq  \frac{4}{3} \text{ \ for all } k \geq 1  \text{ \ \ and \ \ } \srh(k) \leq \frac{2}{\sqrt{3}} \text{ \ for all  } k \geq 20$.
\item[2.]  $\sr(2) = \srh(2)$  and $\sr(3) \geq \srh(3)-0.03$. 
\item[3.]  $\srh(k) \leq \frac{(\log(195 k)+8)^2}{\pi(k-1)}$ \text{ \ for all  } $k \geq 76$. 
\end{itemize}   
\label{thm:intermsys3}
\end{thm}
In fact using the same arguments as in the proof of \textbf{Theorem \ref{thm:intermsys3}}, it can be shown that for $1 \leq k \leq 25$ $\srh(k)$ satisfies the same upper bound as $\sr(k)$ in \textit{Table~\ref{tab:RS_large_sys}}. This leads us to the following conjecture.
\begin{con} Let $\sr(k)$ and $\srh(k)$ be the supremum of the systolic ratio and the homological systolic ratio among all closed, oriented, Riemannian surfaces of genus $k \geq 1$. Then
\[ 
\sr(k) = \srh(k).
\] 
\end{con}
This could in principle be deduced using the same arguments as in the proof of \textbf{Theorem \ref{thm:intermsys3}}. But to this end the upper and lower bound for $\sr(k)$ in any genus would have to be sufficiently close. The idea of Section 2 is to construct new surfaces with large systolic ratio from extremal surfaces. In Section 3 we expand on this idea. If $(M,g)$ and $(N,h)$ are two manifolds of dimension $m$ and $n$, respectively, then
\[
\sy(M \times N,g \times h) = \min( \sy(M,g),\sy(N,h) ).
\]
This enables us to construct manifolds with large systolic ratio from lower dimensional manifolds with large systolic ratio. We illustrate the consequences of this equation by two examples. First we construct $n$-dimensional Euclidean and general product-tori with large systolic ratio from lower dimensional ones, then we construct product manifolds of surfaces and tori. This enables us to prove \textbf{Theorem \ref{thm:lattices}}:

\begin{thm} Let  $\gamma_n$ be Hermite's constant for flat tori in dimension $n$, then
\begin{itemize}
\item[1.]  $\gamma_{m + n} \geq \gamma_{m}^{\left(\frac{m}{m+n}\right)}\cdot \gamma_{n}^{\left(\frac{n}{m+n}\right)}$.
\item[2.]  $\gamma_{2n} \geq \gamma_n$.
\item[3.]  $ \gamma_{n}^{\left(\frac{n}{n-1}\right)} \geq  \gamma_{n+1} \geq \gamma_{n}^{\left(\frac{n}{n+1}\right)}$.
\end{itemize}
\end{thm}
We think this result is known but, as we did not find any proof in the literature (\cite{cs},\cite{ma}), we give one in this paper. Even if the result should be known,  we think that the proof given here illustrates well the techniques used in this paper. More refined methods can be applied to find lower bounds for Hermite's constant for flat tori (see \cite{ma}, p. 92) based on similar ideas. These lead to the laminated lattices, which provide the best known lower bounds for $\gamma_n$ in certain dimensions (see \cite{ma}, \textit{Table 14.4.1}). However, \textbf{Theorem \ref{thm:lattices}} provides practical a priori bounds. Notably, the lower bound in \textbf{Theorem \ref{thm:lattices}-3} completes the known upper bound, which is Mordell's inequality. The same inequalities hold for manifolds homeomorphic to Euclidean tori. These are stated in \textbf{Theorem \ref{thm:homtori}}. Furthermore in \textbf{Theorem \ref{thm:S_T_example}}, we prove:
\begin{thm} Let $b(M,g)=\sum_{i=0}^{m} b_i(M,g)$ be the sum of the Betti numbers of a manifold $(M,g)$ of dimension $m$. Then in each dimension $m \geq 3$ there exist manifolds $R^m_k=(S \times \T^{m-2}, g \times g_E)$, that are product manifolds of a surface $(S,g)$ of genus $k \gg 2^m$ and an Euclidean torus $(\T^{m-2},g_E)$, such that
\[
       C_{1,m} \cdot \frac{\log(b(R^m_k))^2}{b(R^m_k)}   \leq   \sr(R^m_k) \leq C_{2,m} \cdot \frac{\exp(C_{3,m} \sqrt{\log( b(R^m_k))})}{b(R^m_k)}.
\]
\end{thm}

\section*{Acknowledgment}
The second author has been supported by the Alexander von Humboldt foundation. We are grateful to Ivan Babenko, St\'ephane Sabourau and Benjamin Hennion for helpful discussions.

\section{Construction of surfaces with large systolic ratio}

\textbf{proof of Theorem \ref{thm:intermsys}} \\
\\
\textit{1. $\sr(k_1 +k_2) \geq  \left(\frac{1}{\sr(k_1)} + \frac{1}{\sr(k_2)} \right)^{-1}$.}\\ 
\\
Let $\epsilon_1 > 0$ be a positive real number. For $i \in \{1,2\}$, let $(S_i,g_i)$ be a surface of genus $k_i \geq 1$, which satisfies
\begin{equation}
       \sr(S_i,g_i) = \sr(k_i)-\epsilon_1 \text{ \ and \ } \sy(S_i,g_i)=1, \text{ \ hence \ } \vol(S_i,g_i)=\frac{1}{\sr(S_i,g_i)}.
\label{eq:sri}
\end{equation}
To prove our theorem we construct a new surface $(S_c,g_c)$ of genus $k_1 + k_2$ from the surfaces $(S_1,g_1)$ and $(S_2,g_2)$ such that
\[
 \sr(S_c,g_c) = \left(  \frac{1}{\sr(S_1,g_1)}+ \frac{1}{\sr(S_2,g_2)} \right)^{-1}.
\]
As $(S_c,g_c)$ has genus $k_1 + k_2$, we obtain the inequality of \textbf{Theorem \ref{thm:intermsys}-1} from the fact that $\epsilon_1$ can be chosen arbitrarily small.\\
\\
We first construct $(S_c,g_c)$. For fixed $i$ let $a_i$ be a systole of $(S_i,g_i)$. We first divide each $a_i$ into two arcs, $b_i$ and $c_i$ of equal length. Then we cut $a_i$ along $c_i$ and call the surface obtained this way $(S^o_i,g_i)$. We denote by $\alpha_i$ the boundary curve of $(S^o_i,g_i)$. Let $\alpha^1_i$ and $\alpha^2_i$ be the two parts of $\alpha_i$ with common endpoints on $b_i$. We identify the boundary components of  $(S_1^o,g)$ and $(S_2^o,g)$ in the following way
\begin{equation}
    \alpha^{1}_1 \sim \alpha^2_2  { \ \ and \ \ } \alpha^{2}_1 \sim \alpha^1_2 
\label{eq:paste1}
\end{equation}
to obtain a closed surface.

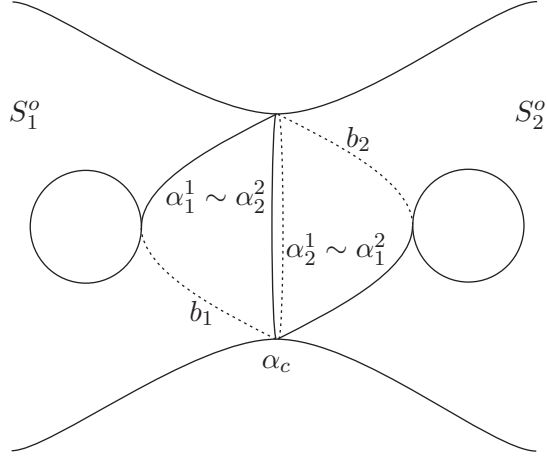
\begin{figure}[h!]
\SetLabels
\L(.28*.70) $S_1^o$\\
\L(.70*.70) $S_2^o$\\
\L(.43*.32) $b_1$\\
\L(.56*.65) $b_2$\\
\L(.41*.54) $\alpha^{1}_1 \sim \alpha^2_2$\\
\L(.51*.44) $\alpha^{1}_2 \sim \alpha^2_1$\\
\L(.49*.23) $\alpha_c$\\
\endSetLabels
\AffixLabels{%
\centerline{%
\includegraphics[height=7cm,width=8cm]{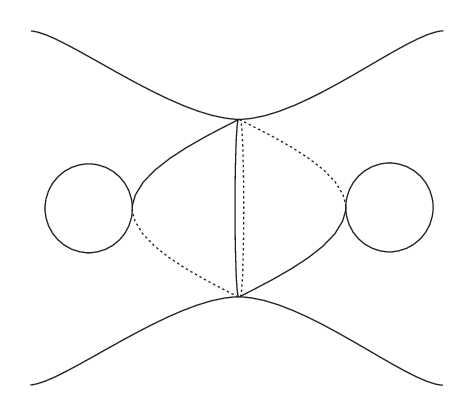}}}
\caption{The surfaces $S_1^o$ and $S_2^o$ with identified boundaries.}
\label{fig:paste}
\end{figure}

We denote the surface of genus $k_1 +k_2$ obtained according to this pasting scheme as
\[
    S_c= S_1^o + S_2^o  \mod (\ref{eq:paste1}). \text{ \ (see Fig.~\ref{fig:paste}) }
\]
We denote by $\alpha_c$ the curve, which is the image of $\alpha_1$ in $S_c$. As the metric on $S_c$, we take the metric of the parts to obtain a surface $(S_c,g_c)$ with the singularity, which is the curve $\alpha_c$.

We now show that the length $\ell(\eta)$ of a non-contractible loop $\eta$ in $(S_c,g_c)$ satisfies
\[
\ell(\eta) \geq \sy(S_1,g_1)= \sy(S_2,g_2)=1.
\]
As it is well-known that every non-contractible loop contains a simple non-contractible sub-loop, i.e. a non-contractible loop without self-intersection, of equal or shorter length, we assume that $\eta$ is a simple closed curve.\\
To prove that $\ell(\eta) \geq 1$ we distinguish two cases: either $\eta$ is contained in either $S_1^o \subset S_c$ or $S_2^o \subset S_c$ or not. Consider the first case.\\
\\
\textit{Case 1: $\eta$ is either contained in $S_1^o$ or contained in $S_2^o$ }\\
\\
We assume without loss of generality that $\eta$ is contained in $S_1^o$. We have to prove that $\ell(\eta) \geq 1$. Now if $\eta$ is non-contractible both in $S_1$ and in $S_1^o$, then $\eta$ is a non-contractible loop in $S_1$ and hence 
\[
    \ell(\eta) \geq  \sy(S_1,g_1)=1
\]
and there is nothing to prove. Therefore it remains to prove the case, where $\eta$ is contractible in $S_1$, but non-contractible in $S_1^o$. It follows from surface topology that if $\delta$ is a closed curve that satisfies this condition, then 
\[
 [\delta] = [\alpha_1]^l \in \pi_1(S_1^o), \text{ \ for some \ } l \in \Z.
\]
As by assumption $\eta$ is additionally a simple loop it follows that $[\eta] = [\alpha_1]^{\pm 1}$. We assume without loss of generality that $[\eta] = [\alpha_1]$. 
We recall that $b_1 \subset S_1^o$ is the part of the systole $a_1 \subset S_1$, that is not cut in $S_1^o$. As $\eta$ runs around the cut in $S_1^o$, whose boundary is $\alpha_1$, it follows that there are two intersection points, $p_1$ and $p_2$ on $b_1$ (see \textit{Fig.~\ref{fig:hole}}), such that
\begin{itemize}
\item[-] $p_1$ and $p_2$ divide $\eta$ into two parts, $\eta_1$ and $\eta_2$, such that 
\item[-] $\eta_1$ is homotopic with fixed endpoints $p_1$ and $p_2$ to an arc $r_1 \alpha^1_1 r_2$, where
\item[-] $r_1$ is the shorter arc on $b_1$ connecting $p_1$ and an endpoint of $\alpha^1_1$ and $r_2$ is  the shorter arc on $b_1$ connecting $p_2$ and an endpoint of $\alpha^1_1$ 
\item[-] $\eta_2$ is homotopic with fixed endpoints to $r_1 \alpha^2_1 r_2$. 
\end{itemize}
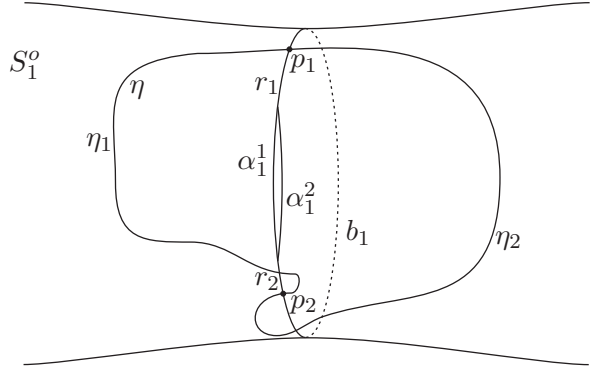
\begin{figure}[h!]
\SetLabels
\L(.25*.80) $S_1^o$\\
\L(.53*.35) $b_1$\\
\L(.35*.74) $\eta$\\
\L(.31*.60) $\,\eta_1$\\
\L(.65*.34) $\,\eta_2$\\
\L(.45*.75) $\,r_1$\\
\L(.45*.23) $\,r_2$\\
\L(.48*.81) $\,p_1$\\
\L(.48*.17) $\,p_2$\\
\L(.44*.54) $\alpha^{1}_1$\\
\L(.48*.44) $\alpha^{2}_1$\\
\endSetLabels
\AffixLabels{%
\centerline{%
\includegraphics[height=5cm,width=8cm]{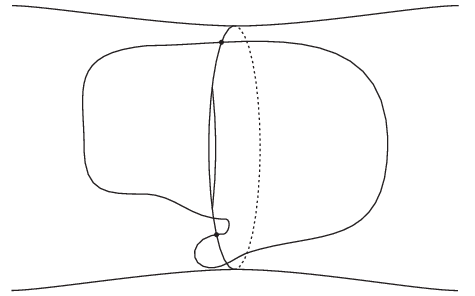}}}
\caption{The surface $S_1^o$ with a curve $\eta$, such that $[\eta] = [\alpha_1]$ in $\pi_1(S_1^o)$.}
\label{fig:hole}
\end{figure}
We now show that 
\begin{equation}
   \ell(\eta_i)  \geq \frac{1}{2} \text{ \ \ for \ \ } i \in \{1,2\},   \text{ \ \ hence \ \ } \ell(\eta) \geq 1. 
\label{eq:eta_S1}
\end{equation}
Let $b'$ be the arc of $b_1$ connecting $p_1$ and $p_2$. We have that $\ell(b') \leq \frac{1}{2}$. 

Furthermore for $i \in \{1,2\}$
\[
\ell(\eta_ib') = \ell(\eta_i) + \ell(b')  \geq \ell(a_1) =1.  
\]
Because otherwise we could in $S_1$ replace $a_1$ by $\eta_ib'$  to obtain a non-contractible loop shorter than $a_1$. But $a_1$ is the systole of $S_1$. A contradiction. Now as $\ell(b') \leq \frac{1}{2}$, it follows that $\ell(\eta_i) \geq \frac{1}{2}$ and hence our statement in (\ref{eq:eta_S1}). This concludes the proof in \textit{Case 1}. \\
\\
\textit{Case 2: $\eta$ is not contained in either $S_1^o$ or $S_2^o$}\\
\\ 
For fixed $i$, we call a loop in $S_c$ \textit{retractable} into $S_i^o \subset S_c$ if and only if it is freely homotopic to a loop contained in $S_i^o$.
We distinguish two subcases: either $\eta$ is retractable into $S_1^o$ or $S_2^o$ or not.\\
\\
\textit{Case 2.a): $\eta$ is retractable into $S_1^o$ or $S_2^o$}\\
\\
Assume without loss of generality that $\eta$ is retractable into $S_1^o \subset S_c$. We first prove the following lemma: 
\begin{lem} For fixed $i \in \{1,2\}$, let $\delta_i \subset S_i^o$ be an arc in $S_c$ with endpoints $p_1$ and $p_2$ on $\alpha_c$, such that $\delta_i$ is homotopic with fixed endpoints to a geodesic arc $d_i \subset \alpha_c$ of length $\ell(d_i) < \ell(\alpha_c) =1$. Then  
\[
    \ell(\delta_i) \geq \ell(d_i).
\]    
\label{thm:comp_bound}
\end{lem}
\textbf{proof of Lemma \ref{thm:comp_bound}} Fix $i$ and let us work in $S_i^o$. It follows from \textit{Case 1} that $\alpha_i$ is a systole of $S_i^o$. Let $b'$ be the remaining arc of $\alpha_i$ connecting $p_1$ and $p_2$, such that $\delta_ib' \subset S_i^o$ is a closed curve. Then $\delta_ib'$ is contained in $S_i^o$ and freely homotopic to the systole $\alpha_i$ of $S_i^o$. It follows that 
\[
  \ell(\delta_i) +  \ell(b') =  \ell(\delta_i b')  \geq  \ell(\alpha_i) = \ell(d_i b') = \ell(d_i) + \ell(b'), \text{ \ \ hence \ \ } \ell(\delta_i) \geq \ell(d_i). 
\]
This proves our lemma. \hfill $\square$\\
Now if $\eta \subset S_c$ is a loop that is retractable into $S_1^o$ then, using \textbf{Lemma \ref{thm:comp_bound}} we can find a comparison curve $\eta'$ for $\eta$, such that 
\[  
     \ell(\eta) \geq \ell(\eta'), [\eta] = [\eta'] \in \pi_1(S_c) \text{ \ \ and \ \ } \eta' \subset S_1^o.
\]
We obtain $\eta'$ from $\eta$ by replacing any arc $v$ of $\eta$ that is contained in either $S_1^o$ or $S_2^o$ and that is homotopic with fixed endpoints to an arc $d \subset \alpha_c$ by the boundary arc $d$.\\ 
That \textbf{Lemma \ref{thm:comp_bound}} can indeed be applied can be seen in the following way: We note that $\eta$ has no self-intersection. Hence $v$ has no self-intersection. Now if the length $\ell(d)$ of $d$ was bigger than one then this would imply that $d$ and hence $v$ has a self-intersection. A contradiction.\\
As $\eta$ is retractable into $S_1^o$ and $\eta'$ is freely homotopic to $\eta$, due to our procedure $\eta'$ is contained in $S_1^o$. By deforming $\eta'$ slightly, we may assume that $\eta'$ is contained in the interior of $S_1^o$. Therefore it follows from \textit{Case 1} that $\ell(\eta') \geq 1$. Hence 
\[
      \ell(\eta) \geq \ell(\eta') \geq 1. 
\]

\textit{Case 2.b): $\eta$ is not retractable into either $S_1^o$ or $S_2^o$}\\
\\
$\eta$ is not retractable into $S_1^o$ and not retractable into $S_2^o$. Now due to this property, $\eta$ contains two subarcs $\eta'_1$ and $\eta'_2$ such that
\begin{itemize}
\item[-] $\eta'_i$ has endpoints on $\alpha_c$ and is contained in $S_i^o$ 
\item[-] there is an arc $b^i$ of $\alpha_c$ connecting the endpoints of $\eta'_i$ on $\alpha_c$ such that  
\item[-] $\eta'_1 b^1$ is not retractable into $S_2^o$ and $\eta'_2 b^2$ is not retractable into $S_1^o$. 
\end{itemize}
It follows from these properties that 
\[
[\eta'_1 b^1] \neq 0 \in \pi_1(S_c) \text{ \ and \ } [\eta'_2 b^2] \neq 0 \in \pi_1(S_c).
\]

We now show that for fixed $i \in \{1,2\}$:
\[
    \ell(\eta'_i) \geq \frac{1}{2}.
\]
Consider without loss of generality $\eta'_1$. Again we distinguish two subcases: the arc $b^1$ is smaller or equal to $\frac{1}{2}$ or not. \\
\\
\textit{Case i): $\ell(b^1) \leq \frac{1}{2}$}\\
\\
As $[\eta'_1b^1] \neq 0 \in \pi_1(S_c)$ and $\eta'_1b^1 \subset S_1^o$, $\eta'_1b^1$ fulfills the conditions of \textit{Case 1} and therefore $\ell(\eta'_1b^1) \geq 1$. As $\ell(b^1) \leq \frac{1}{2}$ it follows that 
\[ 
    \ell(\eta'_1) \geq \frac{1}{2}.
\]
This settles our claim in \textit{Case i}.\\
\\
\textit{Case ii): $\ell(b^1) > \frac{1}{2}$}\\
\\
In this case we close $S_1^o$ along $\alpha_1$ to obtain $S_1$. We denote all curves from $S_1^o$ by the same name in $S_1$. Let $d \subset S_1$ be the shortest geodesic arc on the systole $a_1$ of $S_1$ connecting the endpoints of $b^1$. As  $\ell(b^1) > \frac{1}{2}$, we have that 
\[ 
    \ell(d) < \frac{1}{2}.
\]
Furthermore $[b^1d^{-1}] = 0 \in \pi_1(S^o_1)$. Hence $\eta'_1d$ is a closed curve in $S_1$ that, in $S_1$, is in the same free homotopy class as $\eta'_1b^1$. Now 
\[
[\eta'_1d]=[\eta'_1b^1] \neq 0 \in \pi_1(S_1),
\]
because otherwise it would follow that $[\eta'_1b^1]= [\alpha_1]^{\pm1} \in \pi_1(S_1^o)$ (see \textit{Case 1}). But then $\eta'_1b^1$ would be a curve that is retractable into $S_2^o$, a contradiction.\\ 
Hence $\eta'_1d$ is a non-contractible loop in $S_1$. Its length is bigger or equal to the length of the systole $a_1$ of $S_1$. It follows that
\[
\ell(\eta'_1d) = \ell(\eta'_1) + \ell(d) \geq 1 \text{ \ and \ }  \ell(d) \leq \frac{1}{2}, \text{ \ hence \ } \ell(\eta'_1) \geq \frac{1}{2}.
\]
This settles our claim in \textit{Case ii}.\\
As the same arguments in \textit{Case i} and \textit{Case ii} for $\eta'_1$ apply to $\eta'_2$, we conclude that 
\[ 
    \ell(\eta'_i) \geq \frac{1}{2} \text{ \ \ for \ \ } i \in \{1,2\},  \text{ \ \ hence \ \ } \ell(\eta) \geq \ell(\eta'_1) + \ell(\eta'_2) \geq 1.
\]      
In total we obtain in both \textit{Case 1} and \textit{Case 2} that $\ell(\eta) \geq 1$. \\
As any non-contractible loop in $(S_c,g_c)$  has length greater than or equal to one, we have shown that
\[
\sy(S_c,g_c)=1.  
\]
Due to Equation (\ref{eq:sri}), we have that  $\sr(k_i)-\epsilon_1=\sr(S_i,g_i) = \vol(S_i,g_i)^{-1}$, hence
\[
 \sr(S_c,g_c) = \frac{1}{ \vol(S_1,g_1) + \vol(S_2,g_2)}=\left( \frac{1}{\sr(k_1)-\epsilon_1} + \frac{1}{\sr(k_2)-\epsilon_1 }\right)^{-1}.
\]
Now for every $\epsilon_2 >0$, we can approximate our non-smooth surface $(S_c,g_c)$ with a smooth surface $(S_{\epsilon_2},g_{\epsilon_2})$ such that the distance function and the area of $(S_{\epsilon_2},g_{\epsilon_2})$ is ${\epsilon_2}$-close to that of $(S_c,g_c)$. Letting $\epsilon_1$ and $\epsilon_2$ tend to zero we obtain:
\[
    \sr(k_1 + k_2) \geq \sr(S_c,g_c) \geq \left( \frac{1}{\sr(k_1)} + \frac{1}{\sr(k_2) }\right)^{-1}.
\]
This concludes the proof of \textbf{Theorem \ref{thm:intermsys}-1}. \hfill  $\square$\\
\\
\textit{2. $\srh(k_1 +k_2) \geq  \left(\frac{1}{\srh(k_1)} + \frac{1}{\srh(k_2)} \right)^{-1}$.}\\ 
\\
Let $\epsilon_1 > 0$ be a positive real number. For $i \in \{1,2\}$, let $(S^h_i,g_i)$ be a surface of genus $k_i \geq 1$, which satisfies
\begin{equation}
       \srh(S^h_i,g_i) = \srh(k_i)-\epsilon_1 \text{  \ and \ } \syh(S^h_i,g_i)=1, \text{ \ hence \ } \vol(S^h_i,g_i) =\frac{1}{\srh(S^h_i,g_i)}.
\label{eq:srhi}
\end{equation}
To prove  \textbf{Theorem \ref{thm:intermsys}-2} we construct a new surface $(S^h_c,g_c)$ of genus $k_1 + k_2$ from the surfaces $(S^h_1,g_1)$ and $(S^h_2,g_2)$ such that
\begin{equation}
 \srh(S^h_c,g_c) = \left(  \frac{1}{\srh(S^h_1,g_1)}+ \frac{1}{\srh(S^h_2,g_2)} + \epsilon_2 \right)^{-1} ,
\label{eq:nonsep_srh} 
\end{equation}
where $\epsilon_2$ is a positive real number that can be chosen arbitrarily small.
As $(S^h_c,g_c)$ has genus $k_1 + k_2$, we obtain the inequality of \textbf{Theorem \ref{thm:intermsys}-2} from the fact that $\epsilon_1$ and $\epsilon_2$ can be chosen arbitrarily small.\\
\\
We first construct $(S^h_c,g_c)$. For fixed $i$ let $a_i$ be a homological systole of $(S^h_i,g_i)$.We first divide each $a_i$ into two arcs, $b_i$ and $c_i$, where $c_i$ has length 
\[
\ell(c_i)=\epsilon_2 < \frac{1}{4}.
\]
Then we cut each $a_i$ along the length of $c_i$ and call the surface obtained this way $(S^o_i,g_i)$. We denote by $\gamma_i$ the boundary curve of $(S^o_i,g_i)$.  Now take an Euclidean cylinder $C$, such that 
\begin{itemize}
\item[-] $C$ has height $1$ and $\ell(\partial_1 C) = \ell(\partial_2 C)= \epsilon_2$ , hence $\vol(C)=\epsilon_2$
\item[-] $\gamma$ is the simple closed geodesic in $C$ of length $\ell(\gamma)=\epsilon_2$, which is freely homotopic to the boundary curve $\partial_1 C$ and such that $\dist(\gamma,\partial_1 C) = \frac{1}{2}$.
\end{itemize}
We connect the boundary components of  $(S_1^o,g_1)$ and $(S_2^o,g_2)$ by connecting them with the cylinder $C$ in the following way
\begin{equation}
    \gamma_1 \sim \partial_1 C   { \ \ and \ \ } \gamma_2 \sim \partial_2 C    
\label{eq:paste2}
\end{equation}
to obtain a closed surface.
We denote the surface of genus $k_1 +k_2$ obtained according to this pasting scheme as
\[
    S^h_c= S_1^o +  S_2^o + C  \mod (\ref{eq:paste2}). 
\]
We denote the boundary curves of the embedded cylinder by the same name as in the cylinders itself.\\
We now show that any non-separating loop in $(S^h_c,g_c)$ has length bigger than or equal to one. Let $\eta$  be such a loop. To simplify our proof, we assume that $\eta$ has no self-intersection and that the arcs of $\eta$ contained in $C$ are geodesic, i.e. straight lines. To prove our statement we distinguish two cases: either $\eta$ intersects $\gamma$ transversally or not.\\
\\
\textit{Case 1. $\eta$ intersects $\gamma$ transversally}\\
\\
Due to our assumption, the subarc $\eta' \subset C$ of $\eta$ intersecting $\gamma$ is a straight line. Hence 
$\eta'$ traverses $C$. It follows that its length is bigger than the height of the cylinder $C$, which is  $1$. In this case we have that
\[
    \ell(\eta) \geq \ell(\eta') \geq 1
\]
This settles our proof in \textit{Case 1}.\\
\\
\textit{Case 2. $\eta$ does not intersect $\gamma$ transversally}\\
\\
In the second case, $\eta$ is a non-contractible loop that does not intersect $\gamma$ transversally. Again we distinguish two subcases: $\eta$ is contained in $C$ or not. If $\eta$ is contained in $C$, then $\eta$ is a separating loop, a contradiction to our assumption.\\
If $\eta$ is not contained in $C$ and does not intersect $\gamma$ transversally, then $\eta$ is freely homotopic to a loop $\eta''$, which is contained in the interior of one of the $\left(S_i^o\right)_{i=1,2}$, say $S^o_1$ and such that $\ell(\eta) \geq \ell(\eta'')$. In this case we have that
\[
    \ell(\eta) \geq \ell(\eta'') \geq \syh(S^h_1,g_1) \geq 1.
\]
Here the second inequality follows from the fact that any non-separating simple loop in $S_i^o$ is also a non-separating simple loop in $S^h_i$. This settles our proof in \textit{Case 2}. In total we conclude that $\syh(S^h_c,g_c) \geq  1$.\\
From the homological systolic ratio $\srh(S^h_c,g_c)$ of $(S^h_c,g_c)$ we obtain inequality (\ref{eq:nonsep_srh}). As in (\ref{eq:nonsep_srh}) $\epsilon_1$ and $\epsilon_2$ can be chosen arbitrarily small this yields
\[
    \srh(k_1 + k_2) \geq \srh(S^h_c,g_c) \geq \left( \frac{1}{\srh(k_1)} + \frac{1}{\srh(k_2) }\right)^{-1}.
\]
Here the first inequality follows from the fact that $(S^h_c,g_c)$ can be approximated by a smooth surface. This concludes the proof of \textbf{Theorem \ref{thm:intermsys}}. \hfill $\square$\\
\\
\textbf{proof of Theorem \ref{thm:intermsys2}} The first inequalities in \textbf{Theorem 1.2-1} and \textbf{1.2-2} are a simple consequence of \textbf{Theorem 1.1-1} and \textbf{1.1-2}, respectively. Here we set $k_1=1$ and $k_2=k$ and use the fact that $\sr(1) = \srh(1) = \frac{2}{\sqrt{3}}$.\\ 
We now prove the second inequality in \textbf{Theorem 1.2-1} and then show how to obtain the remaining inequalities in a similar fashion. We have to show that \\
\\
\textit{1. $\frac{1}{\sr(k)} \leq  \frac{1}{\sr(k+1,\hyp)}+ \frac{1}{\pi}$. }\\
\\
Let $\epsilon_1 > 0$ be a positive real number. Let $(S,\hyp)$ be a hyperbolic surface of genus $k+1 \geq 2$, which satisfies
\begin{equation*}
       \sr(S,\hyp) = \sr(k+1, \hyp).
\end{equation*}
As $\sr(k+1,\hyp) = \srh(k+1,\hyp)$, we may assume that $S$ has a systole $\alpha$ which is a non-separating simple closed curve. We can rescale $(S,\hyp)$ to obtain a surface $(S_{max},g_{sc})$ satisfying
\begin{equation}
       \sr(S_{max},g_{sc}) = \sr(k+1, \hyp) \text{ \ and  \ } \sy(S_{max},g_{sc})=1, \text{ \ hence \ } \vol(S_{max},g_{sc})= \frac{1}{\sr(k+1,\hyp)}.
\label{eq:srhom}
\end{equation}
Let $(S^c,g_{sc})$ be the surface which we obtain by cutting open $(S_{max},g_{sc})$ along $\alpha$. As $\alpha$ is non-separating $(S^c,g_{sc})$  has signature $(k,2)$. Let $\alpha_1$ and $\alpha_2$ be the boundary geodesics of $(S^c,g_{sc})$.\\
Let $D$ be a sphere of constant curvature, whose great circles have length $\ell(\alpha)=1$. Let $D_1$ and $D_2$ be the hemispheres which we obtain by cutting $D$ along a great circle. It follows from the geometry of the sphere that
\[
     \vol(D_1)= \vol(D_2) = \frac{\ell(\alpha)^2}{2 \pi} = \frac{1}{2 \pi}.
\]
For fixed $i \in \{1,2\}$, let $\delta \subset D_i$ be a curve connecting two boundary points, $p_1$ and $p_2$ of $D_i$. It follows from the geometry of $D_i$, that there is a comparison boundary arc $\delta'$ of $D_i$, connecting $p_1$ and $p_2$ that is shorter than or of equal length as $\delta$ and such that $\delta (\delta')^{-1}$ is a contractible loop.
\begin{equation}
    \ell(\delta)  \geq \ell(\delta') \text{ \ and \ } [\delta (\delta')^{-1}] = 0 \in \pi_1(D_i), \text{ \ where \ } \delta' \subset \partial D_i.
\label{eq:comp_R}
\end{equation}
To prove our statement, we construct a surface $(S',g')$ of genus $k$ by pasting $D_1$ and $D_2$ along the boundary geodesics of $(S^c,g_{sc})$. Let $\eta$ be a non-contractible simple loop in $(S',g')$. We first show that 
\[
     \ell(\eta) \geq 1, \text{ \ \ hence \ \ }  \sy(S',g') \geq 1.
\]
If $\eta \subset S'$ is contained in $S^c \subset S_{max}$, then $\ell(\eta) \geq \sy(S_{max},g_{sc})=1$ and there is nothing to prove. If $\eta$ is not contained in $S^c$, then $\eta$ intersects $\alpha_1$ or $\alpha_2$ transversally. Then it follows from the comparison statement (\ref{eq:comp_R}) that there is a non-contractible comparison loop $\eta' \subset S^c \subset S_{max}$, such that
\[
\ell(\eta) \geq \ell(\eta').
\]
Hence for any non-contractible simple loop in $S'$ there is a non-contractible loop $\eta'$ in $S_{max}$, whose length is smaller or equal to the length of $\eta$. Therefore 
\[
   \ell(\eta) \geq \ell(\eta') \geq \ell(\alpha) = 1, \text{ \ hence \ } \sy(S',g') \geq 1.
\]
As $\sy(S',g') \geq 1$ and $\vol(S^c,g_{sc})= \sr(k+1,\hyp)^{-1}$ (see  (\ref{eq:srhom})), we obtain for the systolic ratio of $(S',g')$  
\[
   \sr(k) \geq  \sr(S',g') \geq \frac{1}{\vol(S^c,g_{sc}) + \vol(D)}= \left(\frac{1}{\sr(k+1,\hyp)}  + \frac{1}{\pi} \right)^{-1}.
\]
Here the first inequality follows from the fact that $(S',g')$ can be approximated by a smooth surface. The above inequality implies the second inequality in \textbf{Theorem \ref{thm:intermsys2}-1}. This concludes the proof of \textbf{Theorem \ref{thm:intermsys2}-1}. \hfill  $\square$\\
\\
We obtain the second inequality in \textbf{Theorem \ref{thm:intermsys2}-2} by replacing the surface $(S_{max},g_{sc})$ of genus $k+1$ in the previous proof by a surface $(S_{hom},g)$ of genus $k+1$ satisfying 
\[
\srh(k+1) = \srh(S_{hom},g)-\epsilon_1  \text{ \ and  \ } \syh(S_{hom},g)=1.
\]
where $\epsilon_1 > 0$ is a positive real number that can be arbitrarily close to zero. We cut a surface $(S_{hom},g)$ along a homological systole $\beta$. Then we paste two hemispheres of boundary length $\ell(\beta)$ along the boundary curves of the open surface to obtain a surface $(S'_{hom},g_h)$ of genus $k$. Then we apply similar arguments as in the case of the surface $(S',g')$ to obtain
\[
    \syh(S'_{hom},g_h) \geq 1.
\]
The inequality then follows by from a calculation of the homological systolic ratio $\srh(S'_{hom},g_h)$ of $(S'_{hom},g_h)$. \\ 
\\
We obtain inequality in \textbf{Theorem \ref{thm:intermsys2}-3} from the following construction. We cut a surface $(S,g)$ of genus $k+1$, whose systolic ratio is close to the optimal value $\sr(k+1)$ along its homological systole $\beta'$. Then we paste two hemispheres of boundary length $\ell(\beta')$ along the boundary curves of the open surface to obtain a surface $(S'',g'')$ of genus $k$. The inequality then follows from a calculation of the systolic ratio $\sr(S'',g'')$ of $(S'',g'')$.  Here we use the fact that 
\[ 
   \frac{\ell(\beta')^2}{\vol(S,g)} = \frac{\syh(S,g)^2}{\vol(S,g)} \leq \srh(k).
\]
This concludes the proof of \textbf{Theorem \ref{thm:intermsys2}-3} and hence of \textbf{Theorem \ref{thm:intermsys2}} \hfill $\square$\\
\\
\textbf{proof of Theorem \ref{thm:intermsys3}} This proof is very similar to the proof of \cite{gr2}, \textbf{Theorem 2.C}. However our statement is different. We first show how to obtain the first inequality in \textbf{Theorem \ref{thm:intermsys3}-1}. Then we show how to obtain the second inequality and \textbf{Theorem \ref{thm:intermsys3}-2}  by a simple modification of the proof. We show\\
\\
\textit{For all $\epsilon_1 \in (0,10^{-5}], \text{ \ \ \ }  \srh(k)-\epsilon_1 \leq \frac{4}{3}$ for all $k\neq 0$}\\
\\
We prove our statement by induction: As $\sr(1) = \srh(1) = \frac{2}{\sqrt{3}}$, the statement is true for $k=1$. 
We assume that for all $1 \leq k' < k$
\[
   \srh(k')- \epsilon_1 \leq \frac{4}{3} \text{ \ thus \ } \frac{1}{\srh(k')-\epsilon_1} \geq \frac{3}{4}.
\]  
Let $(S_{hom},g)$ be a surface of genus $k > 1$, which satisfies
\begin{equation}
      \vol(S_{hom},g)=1 \text{ \ \ and \ \ } \syh(S_{hom},g)^2 = \srh(k)-\epsilon_1, \text{ \ hence \ } \srh(k)-\epsilon_1  = \syh(S_{hom},g)^2 .
\label{eq:srhom2}
\end{equation}
Let $\alpha$ be a systole of $(S_{hom},g)$. Two cases can occur. Either $\alpha$ is separating or $\alpha$ is non-separating:\\
\\
\textit{Case 1. $\alpha$ is non-separating}\\ 
\\
In this case it follows with Equation (\ref{eq:srhom2}) and as $\syh(S_{hom},g)^2 \geq \sy(S_{hom},g)^2$ that  
\[
    \srh(k)-\epsilon_1   = \syh(S_{hom},g)^2 = \ell(\alpha)^2 = \sy(S_{hom},g)^2 \leq \sr(k) \leq \frac{4}{3}.
\]
\\
\textit{Case 2. $\alpha$ is separating}\\ 
\\
Let $(S^1,g)$ and $(S^2,g)$  be the surfaces of signature $(k_1,1)$ and $(k_2,1)$, which we obtain by cutting open $(S_{hom},g)$ along the systole $\alpha$. Let $\alpha_1$ be the boundary geodesic of $(S^1,g)$ and $\alpha_2$ be the boundary geodesic of $(S^2,g)$. \\
\\
Let $D$ be a sphere of constant curvature, whose great circles have length $\ell(\alpha)$. Let $D_1$ and $D_2$ be the hemispheres which we obtain by cutting $D$ along a great circle. It follows from the geometry of the sphere that
\[
     \vol(D_1)= \vol(D_2) = \frac{\ell(\alpha)^2}{2 \pi}.
\]
To prove our statement, we construct two surfaces $(S^{1p},g_1)$ of genus $k_1$ and $(S^{2p},g_2)$ of genus $k_2$  by pasting $D_1$ and $D_2$ along the boundary geodesics of  $(S^1,g)$ and $(S^2,g)$, respectively. We have

\begin{cla}
For $i \in \{1,2\}$:
\[
     \syh(S^{ip},g_i)^2   \geq \syh(S_{hom},g)^2 =  \srh(k)-\epsilon_1.
\]
\label{thm:div_Shom}
\end{cla}
\begin{proof} Consider without loss of generality the surface $(S^{1p},g_1)$. Let $\beta$ be a non-separating loop such that $\ell(\beta) = \syh(S^{1p},g_1)$. If $\beta$ is contained in $(S^1,g) \subset (S^{1p},g_1)$ then $\beta$ is also a non-separating loop in $(S_{hom},g)$ and 
\[
       \ell(\beta) = \syh(S^{1p},g_1) \geq  \syh(S_{hom},g).
\]
Hence our statement is true. If $\beta$ is not contained in $(S^1,g) \subset (S^{1p},g_1)$ then some of its arcs traverse the hemisphere $D_1$. Let $\delta \subset D_1$ be a curve connecting two boundary points, $p_1$ and $p_2$ of $D_1$. It follows from the geometry of $D_1$, that there is a comparison boundary arc $\delta'$ of $D_1$, connecting $p_1$ and $p_2$ that is shorter than or of equal length as $\delta$ and such that $\delta (\delta')^{-1}$ is a contractible loop. Hence there is a comparison curve $\beta'$ for $\beta$ in the same homology class as $\beta$ that is contained in $(S^1,g) \subset (S^{1p},g_1)$ of smaller or equal length. Again we conclude that 
\[
         \ell(\beta) \geq \ell(\beta') \geq \syh(S^{1p},g_1) \geq  \syh(S_{hom},g).
\] 
The same arguments for $\syh(S^{2p},g_2)$ yield our claim. 
\end{proof}
\color{black} It follows from \textbf{Claim \ref{thm:div_Shom}} that for $i \in \{1,2\}$:
\[
   \srh(k_i) \geq \srh(S^{ip},g_i) = \frac{  \syh(S_{hom},g)^2 }{\vol(S^{ip},g_i) } \text{ \ thus \ }   \frac{1}{\srh(k_i)}  \leq \frac{\vol(S^{ip},g_i) }{ \syh(S_{hom},g)^2  }.
\]
Combining the above two inequalities and using the fact that $\vol(S^{1p},g_1)+ \vol(S^{2p},g_2)= 1 + \frac{\ell(\alpha)^2}{\pi}$, we have that 
\[
    \frac{1}{\srh(k_1)} + \frac{1}{\srh(k_2)}  \leq  \frac{1+ \frac{\ell(\alpha)^2}{ \pi}}{\srh(k) -\epsilon_1 } 
\]
Now  $\ell(\alpha)^2 = \sy(S_{hom},g)^2 = \sr(S_{hom},g) \leq \sr(k)$. This yields
\begin{equation} 
   \frac{1}{\srh(k_1)} + \frac{1}{\srh(k_2)}   \leq \frac{1+ \frac{\sr(k)}{\pi}}{\srh(k)-\epsilon_1}. 
\label{eq:sep_alpha}   
\end{equation}
Applying the induction hypothesis and using the fact that $\sr(k) \leq \frac{4}{3}$, we obtain:
\[
   \frac{3}{2} \leq  \frac{1}{\srh(k_1)} + \frac{1}{\srh(k_2)}   \leq \frac{1+ \frac{4}{3\pi}}{\srh(k)-10^{-5}}, \text{ \ hence \ }  0.96 \geq \srh(k). 
\]
But this proves our hypothesis. This settles our claim in \textit{Case 2} and therefore concludes the proof of the first part of \textbf{Theorem \ref{thm:intermsys3}-1}. Letting $\epsilon_1$ go to zero, we obtain the first part of \textbf{Theorem \ref{thm:intermsys3}}.  \hfill $\square$ \\
To prove the second part, we use the same arguments. Here we use the fact that we already know that $\sr(k) \leq \frac{2}{\sqrt{3}}$ for $k \geq 20$. In fact using the same arguments, it can be shown that $\srh(k)$ satisfies the same upper bound as $\sr(k)$ in \textit{Table~\ref{tab:RS_large_sys}}. \\
To prove the first part of \textbf{Theorem \ref{thm:intermsys3}-2}, we also follow the above proof. However in this case, we obtain a contradiction in  inequality (\ref{eq:sep_alpha}) of \textit{Case 2} from the known value of $\srh(1)$ and the upper and lower bound for $\sr(2)$ and $\srh(2)$ (see \textit{Table \ref{tab:RS_large_sys}}). Hence this case leads to a contradiction. It remains \textit{Case 1} from which follows that  
\[
\sr(2) = \srh(2).
\]
We prove the second part,
\[
\sr(3) \geq \srh(3) -0.03
\]
in a similar fashion. If \textit{Case 2} holds, then we use iteratively the second inequality in (\ref{eq:sep_alpha}) to show that in this case $\srh(3) \leq 0.69$, from which follows our statement. This settles the proof of the first and second part of \textbf{Theorem 1.3} \hfill $\square$ \\
\\
\textbf{proof of Theorem 1.3-3} To prove the third part of the theorem we first state a good upper bound for $\sr(k)$ which is proven in the appendix: 
\begin{equation} 
    \sr^u(k) := 
\left\{ {\begin{array}{*{20}c}
   {4/3}  \\
   {\frac{(\log(50 k)+1.4)^2}{\pi(k-1)} }  \\
\end{array}} \right.\text{ \ if \ }\begin{array}{*{20}c} 
   k < 17  \\
   k \geq 17 \\
\end{array}  \text{ \ we have \ } \sr^u(k) \geq \sr(k). \text{ \ \ (see appendix)}
\label{eq:sruk}
\end{equation}

Set, for $k \in [0,\infty)$
\[ 
    \sr^u_h(k) := 
\left\{ {\begin{array}{*{20}c}
   {4/3}  \\
   {\frac{(\log(195 k)+8)^2}{\pi(k-1)} }  \\
\end{array}} \right.\text{ \ if \ }\begin{array}{*{20}c}
   k \leq 75  \\
   k > 75. \\
\end{array}
\]
We now prove by induction that
\[
\srh(k) \leq \sr^u_h(k) \text{ \ \ for all \ \ } k \in \N \backslash \{0\}.
\]
Therefore we use the same arguments as in the proof of the first part of \textbf{Theorem \ref{thm:intermsys3}}. For fixed $k \in \N \backslash \{0\}$, we assume that our statement is proven for all $k' <  k$. It is easy to see that the crucial point is inequality (\ref{eq:sep_alpha}) in \textit{Case 2}, which states that
\[
   \frac{1}{\srh(k_1)} + \frac{1}{\srh(k_2)}   \leq \frac{1+ \frac{\sr(k)}{\pi}}{\srh(k)}, \text{ \ where \ } k_1+k_2 =k.
\]
Applying the induction hypothesis and the upper bound for $\sr(\cdot)$ and $\srh(\cdot)$ , we obtain from the above inequality:
\begin{equation}
   \frac{1}{\sr^u_h(k_1)} + \frac{1}{\sr^u_h(k_2)}  \leq  \frac{1}{\srh(k_1)} + \frac{1}{\srh(k_2)}   \leq \frac{1+ \frac{\sr^u(k)}{\pi}}{\srh(k)}. 
\label{eq:sep_alpha_log}   
\end{equation}
It remains to show that 
\begin{equation}
\frac{1}{\sr^u_h(k)}  \leq  \frac{\frac{1}{\sr^u_h(k_1)} + \frac{1}{\sr^u_h(k_2)}}{1+ \frac{\sr^u(k)}{\pi}} \text{ \ \ (see appendix).}
\label{eq:sruh_srh}
\end{equation}
This implies $\sr^u_h(k) \geq \srh(k)$ and hence our hypothesis is true.\\ 
Inequality (\ref{eq:sruh_srh}) is shown in the appendix. This settles our claim in \textit{Case 2} and therefore concludes the proof of \textbf{Theorem \ref{thm:intermsys3}-3}. In total we have proven \textbf{Theorem \ref{thm:intermsys3}}.  \hfill $\square$

\section{Construction of manifolds with large systolic ratio}

In this section we construct manifold with large systolic ratio from lower dimensional manifolds with large systolic ratio. To this end we first prove the following lemma:

\begin{lem} Let $(M,g)$ and $(N,h)$ be two closed, oriented Riemannian manifolds of dimension $m$ and $n$, respectively. If the product manifold $(M \times N,g \times h)$ has a systole, we have:
$$\sy(M \times N,g \times h) = \min( \sy(M,g),\sy(N,h) ).$$
\label{thm:minNM}
\end{lem}
\textbf{proof of Lemma \ref{thm:minNM}}
Let $\eta$ be a systole of $(M \times N,g \times h)$ of length
\[
\ell(\eta)=\sy(M \times N,g \times h).
\]
Let $p_M$ an $p_N$ be the canonical projection from $(M \times N,g \times h)$ to $(M,g)$ and $(N,h)$, respectively. Consider the curves $p_M(\eta)$ and $p_N(\eta)$, respectively. We have:
\[
    \ell( \eta) \geq \ell( p_M(\eta)) \text{ \ \ and \ \ }  \ell( \eta) \geq \ell( p_N(\eta)).
\]
Let $[\eta] \in \pi_1(M \times N)$ be the homotopy class of $\eta$. As $\eta$ is non-contractible we have that $[\eta] \neq 0$. As
\[
\pi_1(M \times N)= \pi_1(M) \times \pi_1(N) = {p_M}_*(\pi_1(M \times N)) \times {p_N}_*(\pi_1(M \times N)),
\]
either $[p_M(\eta)] \neq 0$ or $[p_M(\eta)] \neq 0$. Assume without loss of generality that $[p_M(\eta)] \neq 0$. It follows that
\[
     \sy(M \times N,g \times h) = \ell( \eta) \geq \ell( p_M(\eta)) \geq \sy(M,g) \geq \min( \sy(M,g),\sy(N,h) ).
\]
This proves our lemma. \hfill $\square$\\
\\
As a corollary we obtain:
\begin{cor} Let $(M,g)$ and $(N,h)$ be two closed, oriented Riemannian manifolds of dimension $m$ and $n$, respectively, such that $\sy(M,g)=\sy(N,h)$. Then for the product manifold $(M \times N,g \times h)$, we have:
\[
   \sr(M \times N,g \times h) = \sr(M,g)\cdot\sr(N,h).
\]
\label{thm:sr_prod}
\end{cor}
\textbf{proof of Corollary \ref{thm:sr_prod}} This follows immediately from \textbf{Lemma \ref{thm:minNM}} and the fact that \\
$\vol(M \times N,g \times h)=\vol(M,g)\cdot\vol(N,h)$. \hfill $\square$\\
\\
Note that we can always scale one of the two manifolds in the product manifold to meet the conditions of the corollary. In the following we apply the corollary to two examples. First we construct $n$-dimensional Euclidean product-tori with large systolic ratio from lower dimensional ones, then we construct product manifolds of surfaces and tori. \\
\\
\textbf{Example 1 - tori} A \textit{lattice} $\Lambda$ of dimension $n$ is a discrete subgroup of $\R^n$ that spans $\R^n$. An $n$-dimensional \textit{flat torus} $\T^n = \R^n / \Lambda$ is the quotient of $\R^n$ and a lattice $\Lambda$.
The shortest non-zero lattice vector of $\Lambda$ is the systole of $\T^n$. It's length $\sy(\T^n)$ is
\[
    \sy(\T^n) = \min \limits_{\lambda \in \Lambda \backslash \{0\}} \| \lambda \|_2,
\]
where $\| \cdot \|_2$ denotes the Euclidean norm. If $A$ is a matrix representation of a basis of $\Lambda$, then $\det(\Lambda)$, the \textit{determinant of $\Lambda$} is equal to $\det(A)$ and
\[
\vol(\T^n) = |\det(\Lambda)| = |\det(A)|.
\]
\textit{Hermite's constant or invariant} $\gamma_n$ is given by
\[
    \gamma_n = \max \{ \sr(\T^n)^{\frac{2}{n}} \mid \T^n \text{ \ flat torus of dimension } n \}
\]
It follows from this definition that it is the maximal value that the squared norm of the shortest non-zero lattice vector can attain among all lattices of determinant $1$. Let $B^n$ be a $n$-dimensional Euclidean ball of radius $1$. It was proven by Hlawka \cite{hl} and Minkowski \cite{mi} that
\begin{equation}
\frac{n}{2 \pi e} \leq \left(\frac{\vol(B^n)}{2} \right)^{-\frac{2}{n} } \leq  \gamma_{n} \leq  4\cdot \left(\vol(B^n) \right)^{-\frac{2}{n} } \leq \frac{2n}{\pi e}.
\label{eq:Hermite}
\end{equation}
Here the approximations of the bounds apply for large $n$. From \textbf{Corollary \ref{thm:sr_prod}}, we obtain:
\begin{thm} Let  $\gamma_n$ be Hermite's constant for flat tori in dimension $n$, then
\begin{itemize}
\item[1.]  $\gamma_{m + n} \geq \gamma_{m}^{\left(\frac{m}{m+n}\right)}\cdot \gamma_{n}^{\left(\frac{n}{m+n}\right)}$.
\item[2.]  $\gamma_{2n} \geq \gamma_n$.
\item[3.]  $ \gamma_{n}^{\left(\frac{n}{n-1}\right)} \geq  \gamma_{n+1} \geq \gamma_{n}^{\left(\frac{n}{n+1}\right)}$.
\label{thm:lattices}
\end{itemize}
\end{thm}

Let furthermore $(\T^n,g)$ be a manifold homeomorphic to an Euclidean torus of dimension $n$ or shortly a \textit{torus}. We define \textit{Hermite's constant for general tori} $\delta_n$ by
\[
   \delta_n = \sup \{ \sr(\T^n,g)^{\frac{2}{n}} \mid (\T^n,g) \text{ \ torus of dimension } n \}. 
\]
As $\gamma_n \leq \delta_n$ the same lower bound as in inequality (\ref{eq:Hermite}) applies. An upper bound of order $n$ was conjectured in \cite{gr3}. The best known upper bound is of order $n^2$ and is stated in \cite{na},\textbf{Theorem 4.2}. This implies for large $n$ that
\begin{equation}
     \frac{n}{2 \pi e}   \leq \delta_{n}  \leq \left(\frac{n}{e}\right)^2.
\label{eq:Herm_general}
\end{equation}
Using the same methods as in the Euclidean case we show that 
\begin{thm} Let $\delta_n$ be Hermite's constant for general tori of dimension $n$, then
\begin{itemize}
\item[1.]  $\delta_{m + n} \geq \delta_{m}^{\left(\frac{m}{m+n}\right)}\cdot \delta_{n}^{\left(\frac{n}{m+n}\right)}$.
\item[2.]  $\delta_{2n} \geq \delta_n$.
\label{thm:homtori}
\end{itemize}
\end{thm}

The idea of the proof is to construct tori with large systoles from lower dimensional extremal ones.\\
\\
\textbf{proof of Theorem \ref{thm:lattices}}. Let $\T^{m}$ and $\T^{n}$ be two flat tori of dimension $m$ and $n$, respectively, such that
\[
     \sr(\T^{m})^{\frac{2}{m}} = \gamma_m \text{ \ \ and \ \ }   \sr(\T^{n})^{\frac{2}{n}} = \gamma_n.
\]
We obtain \textbf{Theorem \ref{thm:lattices}-1} by scaling $\T^{m}$ to obtain a torus $\T_1^{m}$ and $\T^{n}$ to obtain a torus $\T_2^{n}$, such that
\[
      \sy(\T_1^{m})=  \sy(\T_2^{n})= 1.
\]
Let $\T_1^{m} \times \T_2^{n}$ be the product torus of dimension $m+n$. Applying \textbf{Corollary \ref{thm:sr_prod}} to $\T_1^{m} \times \T_2^{n}$ we obtain
\[
    {\gamma_{m+n}}^\frac{m+n}{2} \geq \sr(\T_1^{m} \times \T_2^{n}) =  \sr(\T_1^{m})\cdot\sr(\T_2^{n})= \sr(\T^{m})\cdot\sr(\T^{n})= \gamma_m^{\frac{m}{2}}\cdot \gamma_n^{\frac{n}{2}}.
\]
This inequality is equivalent to \textbf{Theorem \ref{thm:lattices}-1}. The second inequality of the theorem follows by setting $n=m$ in the first inequality. The first inequality in \textbf{Theorem \ref{thm:lattices}-3} is Mordell's inequality (see \cite{ma} or \cite{mo}). The second inequality follows by setting $m=1$ in \textbf{Theorem \ref{thm:lattices}-1} and using the fact that $\gamma_1 = 1$. \\
This concludes the proof of \textbf{Theorem \ref{thm:lattices}}. \hfill $\square$\\
\textbf{Theorem \ref{thm:homtori}} follows by the same arguments replacing $\T^n$ by  $(\T^n,g)$ and will not be shown here.\\
\\
\textbf{Example 2 - products of tori and surfaces} Denote from here on by $g_E$ be the Euclidean metric tensor. In this example we construct product manifolds with large systolic ratio from surfaces and tori with large systolic ratio. This enables us to prove:
\begin{thm} Let $b(M,g)=\sum_{i=0}^{m} b_i(M,g)$ be the sum of the Betti numbers of a manifold $(M,g)$ of dimension $m$. Then in each dimension $m \geq 3$ there exist manifolds $R^m_k=(S \times \T^{m-2}, g \times g_E)$, that are product manifolds of a surface $(S,g)$ of genus $k \gg 2^m$ and an Euclidean torus $(\T^{m-2},g_E)$, such that
\[
       C_{1,m} \cdot \frac{\log(b(R^m_k))^2}{b(R^m_k)}   \leq   \sr(R^m_k) \leq C_{2,m} \cdot \frac{\exp(C_{3,m} \sqrt{\log( b(R^m_k))})}{b(R^m_k)}.
\]
\label{thm:S_T_example}
\end{thm}
Here the upper bound is the universal upper bound stated in \cite{sa3}, \textbf{Theorem 1.2}, inequality (1.5).\\
\\
\textbf{proof of Theorem \ref{thm:S_T_example}} Let $(M,g)$ be an $m$-dimensional manifold and let
\[
    P_M(x):= \sum_{i=0}^{m} b_i(M,g)\cdot x^i
\]    
be the \textit{Poincar\'e polynomial}. It is well-known that if  $S=(S,g)$ is a surface of genus $k$ and $T=(\T^{m-2},g_E)$ is an Euclidean torus of dimension $m-2$, then
\[
        P_S(x) = 1 + 2k \cdot x + x^2     \text{ \ \ and \ \ }  P_T(x) = \sum_{i=0}^{m-2}  \binom {m-2} {i} x^i, \text{ \ \ where \ \ }  \sum_{i=0}^{m-2}  \binom {m-2}{i} = 2^{m-2}.
\]
It follows furthermore from the K\"unneth theorem, that for $S\times T =(S \times \T^{m-2}, g \times g_E)$, we have
\[
     P_{S \times T}(x) =  P_S(x)\cdot  P_T(x).
\]
It is easy to deduce from this formula, that 
\begin{equation}
b(S \times T)= 4k + o(1), \text{ \ \ if \ \ } k \gg 2^m.
\label{eq:brmk}
\end{equation}
We now choose our product manifold $R^m_k=(S \times \T^{m-2}, g \times g_E)$ in the following way. Let $\epsilon > 0$ be a positive real number and let in $R^m_k$, $(S,g)$ be a surface of genus $k \gg 2^m$, such that
\[
        \sy(S,g) = 1 \text{ \ \ and \ \ }  \sr(S,g) = \sr(k) - \epsilon.
\]
In $R^m_k$ let furthermore $\T^{m-2}$ be a Euclidean torus of dimension $m-2$, such that
\[
    \sy(\T^{m-2}) = 1 \text{ \ \ and \ \ }  \sr(\T^{m-2}) = \gamma_{m-2}^{\frac{{m-2}}{2}}.
\]
It follows from \textbf{Corollary \ref{thm:sr_prod}} and the inequalities (\ref{eq:srk_bound}) and (\ref{eq:Hermite}) that
\[
      \sr(S \times \T^{m-2}, g \times g_E) = \sr(S, g) \cdot \sr(\T^{m-2},g_E)  \geq K\frac{\log(k)^2}{k} \cdot {(m-2)}^{\frac{{m-2}}{2}}.
\]
Then the lower bound in \textbf{Theorem \ref{thm:S_T_example}} follows by applying Equation (\ref{eq:brmk}) to the above inequality. \hfill $\square$\\

\section{Appendix}
In this part we show the missing inequalities from the proof of \textbf{Theorem \ref{thm:intermsys3}-3}.\\
\\
\textit{ Inequality (\ref{eq:sruk})}\\
\\
We first recall inequality (\ref{eq:srk_detail}): for all  $r \in (0,\frac{1}{8}),\sr(k) \in (0,\frac{4}{3}]$  we know that 
\[
f(r,\sr(k)):= \frac{\log(2r^2 \sr(k))^2}{4\pi\sr(k)\left(\frac{1}{2}-4r\right)^2(k-1)} \geq 1,
\]
where the bound on $\sr(k)$ is the well-known $\frac{4}{3}$ bound. 
Set (for $k \in [0,\infty)$)
\[ 
    \sr^u(k) := 
\left\{ {\begin{array}{*{20}c}
   {4/3}  \\
   {\frac{(\log(50 k)+1.4)^2}{\pi(k-1)} }  \\
\end{array}} \right.\text{ \ for \ }\begin{array}{*{20}c}
   k < 17  \\
   k \geq 17. \\
\end{array}
\]
Let $k \geq 17$. That $\sr^u(k) \geq \sr(k)$ can be seen in the following way: 
\begin{itemize}
\item[1.)] Fixing $r$ and deriving the function $f(r,\cdot)$ we see that this is a monotonically decreasing function in the interval $(0,\frac{1}{2r^2}]$. As $r \in (0,\frac{1}{8})$ this implies that $f(r,\cdot)$ is a monotonically decreasing function in the interval $(0,\frac{4}{3}]$. 
\item[2.)] Examining $f(r,y)$ in the interval $(0,\frac{1}{8}) \times (0,\frac{4}{3}]$ we obtain a value for $r$ that gives a good estimate on $\sr(k)$:  we choose
\[
      r'=\frac{1}{4(\log(25k-25)+2)} < 1/8
\]      
and let $\sr(r',k)$ be the value given by $f(r',\sr(r',k))=1$. Then it follows from (\ref{eq:srk_detail}) that 
\[
         \sr(r',k) \geq \sr(k).
\]         
\item[3.)] It follows from elementary, but tedious, calculations that for $k \geq 17$, 
\begin{equation}
f(r',\sr^u(k)) < 1 \text{ \ \ and \ \ } \lim_{k \to \infty} f(r',\sr^u(k)) = 1 \text{ \ \ \  (see Fig. \ref{fig:frprime})}.
\label{eq:goodsruk}
\end{equation}
A plot of $f(r',\sr^u(k))$ is shown below:
\begin{figure}[h!]
\SetLabels
\L(.82*.08) $k$\\
\L(.82*.62) $f(r',\sr^u(k))$\\
\endSetLabels
\AffixLabels{%
\centerline{%
\includegraphics[height=6cm,width=10cm]{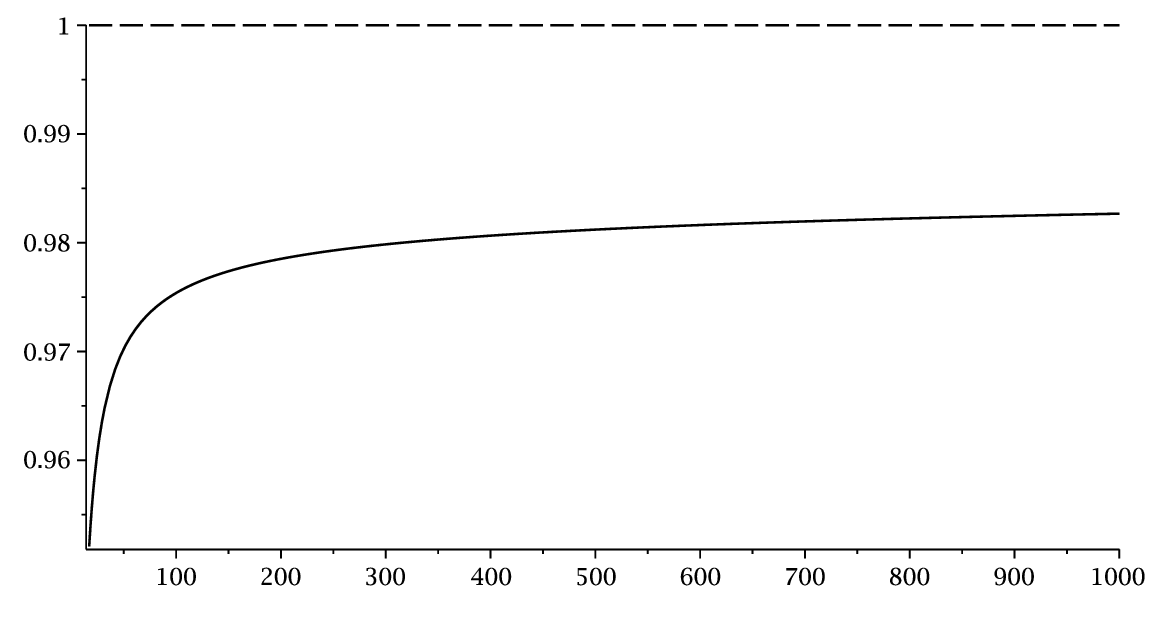}}}
\caption{Plot of the function $f(r',\sr^u(k))$ in the interval $ k \in [17,1000]$.}
\label{fig:frprime}
\end{figure}

But this implies that $\sr^u(k) \geq \sr(r',k) \geq \sr(k)$ for all $k \geq 17$. 
\end{itemize}

\begin{figure}[h!]
\SetLabels
\L(.65*.04) $k$\\
\L(.33*.37) $k_1$\\
\endSetLabels
\AffixLabels{%
\centerline{%
\includegraphics[height=7cm,width=10cm]{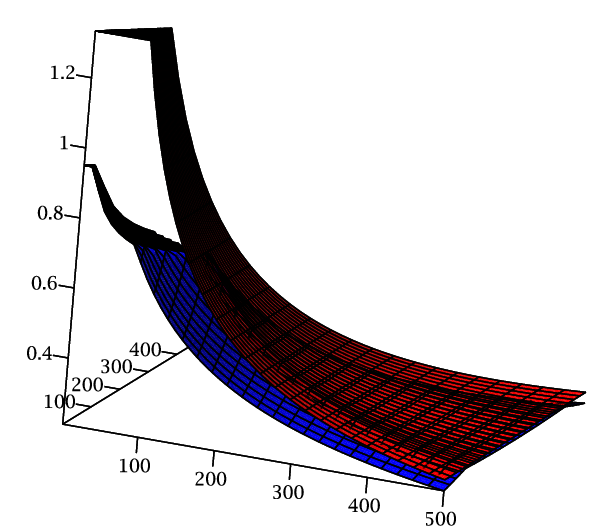}}}
\caption{Plot of the functions $\sr^u_h(k)$ (red) and $s(k,k_1)$ (blue)  in the set $\{ (k,k_1) \in \R^2 | k \in [2,500], k_1 \in [1,k-1]$\}.}
\label{fig:skk1}
\end{figure}

\textit{Inequality (\ref{eq:sruh_srh})} \\
\\
We recall that for $k \in [0,\infty)$
\[ 
    \sr^u_h(k) := 
\left\{ {\begin{array}{*{20}c}
   {4/3}  \\
   {\frac{(\log(195 k)+8)^2}{\pi(k-1)} }  \\
\end{array}} \right.\text{ \ if \ }\begin{array}{*{20}c}
   k \leq 75  \\
   k > 75. \\
\end{array}
\]
It remains to show that for all $k \in \N \backslash \{0,1\}$ and  $1 \leq k_1 \leq k-1$

\begin{equation}
\sr^u_h(k) \geq  \frac{1+ \frac{\sr^u(k)}{\pi}}{\frac{1}{\sr^u_h(k_1)} + \frac{1}{\sr^u_h(k-k_1)}} := s(k,k_1), 
\end{equation}
which is equal to inequality (\ref{eq:sruh_srh}).\\
\\
\textit{Case 1: $k \leq 500$}\\
\\
For $k \leq 500$  it can be verified by calculating $\sr^u_h(k_1)$ and $\sr^u_h(k - k_1)$ explicitly that the inequality is true. A plot of $\sr^u_h(k)$ and $s(k,k_1)$ is shown in Fig. \ref{fig:skk1}.\\

\textit{Case 2: $k > 500$}\\
\\
We now look at a fixed $k > 500$. In this case $1+ \frac{\sr^u(k)}{\pi}$ is constant. To find the maximum of $s(k,k_1)$ we have to find the minimum of the function  
\[
u_k(k_1):=\frac{1}{\sr^u_h(k_1)} + \frac{1}{\sr^u_h(k-k_1)} \text{ \ for \ } k_1 \in [1,k-1].
\]
We now show that $u_k(k_1)$ is minimal for $k_1=75$ and $k_1=k-75$. As $u_k(k_1) = u_k(k-k_1)$ it is sufficient to examine the function in the interval $[1,\frac{k}{2}]$. Furthermore, by definition, 
\begin{itemize}
\item[1.)] $\frac{1}{\sr^u_h(x)}$ is a monotonically increasing function for $x \in [0,\infty)$.
\item[2.)] Hence $\frac{1}{\sr^u_h(k-x)}$ is a monotonically decreasing function for $x \in [0,k]$. 
\item[3.)] As $\frac{1}{\sr^u_h(x)}$ is equal to $\frac{3}{4}$ in the interval $x \in [0,75]$, this implies that  $u_k(x)$ has a local minimum at $x = 75$.
\item[4.)] For $x > 75$ we have that $\frac{1}{\sr^u_h(x)} =  \frac{\pi(x-1)}{\log(e^8 \cdot 195 x)^2} := w(x)$. 
\item[5.)] It can be shown that $w''(x) < 0$ for $x \in [0,\infty)$. This implies that $w(x)$ is a concave function. It follows that $u_k''(x) = w''(x) + w''(k-x) < 0$ is a concave function in the interval $[75,k-75]$.  
\end{itemize}
From 3.) and 5.) we conclude that the local minima of $u_k(x)$ in $x=75$ and $x=k-75$ are indeed global minima in the interval $[0,k]$. It follows that for fixed $k$ and all $k_1 \in [1,k-1]$ 
\[
\frac{1+ \frac{\sr^u(k)}{\pi}}{\frac{3}{4} + \frac{1}{\sr^u_h(k-75)}} \geq \frac{1+ \frac{\sr^u(k)}{\pi}}{\frac{1}{\sr^u_h(k_1)} + \frac{1}{\sr^u_h(k-k_1)}}. 
\]

It remains to show that 
\[
     \sr^u_h(k)  \geq  \frac{1+ \frac{\sr^u(k)}{\pi}}{\frac{3}{4} + \frac{1}{\sr^u_h(k-75)}}.
\]
For $k \in [2,500]$ it can be seen in Fig. \ref{fig:skk1} that the inequality is true.\\
For $k \geq 500$ we plug in the corresponding formulas for $\sr^u(k),\sr^u_h(k)$ and $\sr^u_h(k-75)$ this simplifies to 
\[
  \sr^u_h(k) =  \frac{ \ln(e^{8} \cdot 195  k)^2 }{\pi(k-1)} \geq   \frac{ \ln(e^{8} \cdot 195  (k-75))^2 }{\pi(k-1)} \cdot \left( \frac{ \pi^2(k-1) + \ln(e^{1.4} \cdot 50 k)^2}{\pi ( \pi(k-76) +\frac{3}{4} \cdot \ln(e^{8} \cdot 195  (k-75))^2 )} \right).
\]  
As $\frac{ \ln(e^{8} \cdot 195  k)^2 }{\pi(k-1)} >  \frac{ \ln(e^{8} \cdot 195  (k-75))^2 }{\pi(k-1)}$ we are done when we have shown that 
\[
1 \geq  \frac{ \pi^2(k-1) + \ln(e^{1.4} \cdot 50 k)^2}{\pi ( \pi(k-76) +\frac{3}{4} \cdot \ln(e^{8} \cdot 195  (k-75))^2 )} 
\]
This is  equal to 
\[
  \frac{3 \cdot \pi}{4} \cdot \ln(e^{8} \cdot 195  (k-75))^2 -\pi^2 \cdot 75 \geq \ln(e^{1.4} \cdot 50 k)^2.
\]  
Finally, we note that for $k =500$ the above inequality is true. Deriving both sides of the above inequality we obtain: 
\[
          \frac{3 \cdot \pi}{4} \cdot \ln(e^{8} \cdot 195  (k-75))^2 -\pi^2 \cdot 75)' \geq (\ln(e^{1.4} \cdot 50 k)^2)'.
\]
This inequality is also true for all $k \geq 500$. These two conditions imply that
 inequality (\ref{eq:sruh_srh}) holds for $k \geq 500$. In total we have proven inequality (\ref{eq:sruh_srh}).

\vspace{2cm}

\noindent Hugo Akrout\\	
\noindent Department of Mathematics, Universit\'e Montpellier 2 \\
\noindent place Eug\`ene Bataillon, 34095 Montpellier cedex 5, France \\
\noindent e-mail: \textit{akrout@math.univ-montp2.fr}\\
\\
\noindent Bjoern Muetzel\\	
\noindent Department of Mathematics, Dartmouth College\\
\noindent 27 N. Main street, Hanover, NH 03755, USA \\
\noindent e-mail: \textit{bjorn.mutzel@gmail.com}

\end{document}